\documentclass[a4paper,10pt]{amsart}
\usepackage{amsmath,amssymb,amsxtra,amsfonts,amsthm,url,comment,microtype}

\usepackage{a4wide}

\usepackage{tikz} 
\usetikzlibrary{arrows,backgrounds,decorations,automata}

\theoremstyle{plain}
\newtheorem{theorem}{Theorem}
\newtheorem{lemma}{Lemma}
\newtheorem{corollary}{Corollary}
\newtheorem{proposition}{Proposition}

\theoremstyle{definition}
\newtheorem{remark}{Remark}

\newcommand{\res}[3]{\underset{{{#1}={#2}}}{\textnormal{Res}}\left\{#3\right\}}

\makeatletter
\DeclareRobustCommand\bigop[2][1]{%
  \mathop{\vphantom{\sum}\mathpalette\bigop@{{#1}{#2}}}\slimits@
}
\newcommand{\bigop@}[2]{\bigop@@#1#2}
\newcommand{\bigop@@}[3]{%
  \vcenter{%
    \sbox\z@{$#1\sum$}%
    \hbox{\resizebox{\ifx#1\displaystyle#2\fi\dimexpr\ht\z@+\dp\z@}{!}{$\m@th#3$}}%
  }%
}
\makeatother

\makeatletter
  \def\vhrulefill#1{\leavevmode\leaders\hrule\@height#1\hfill \kern\z@}
\makeatother

\newcounter{nootje}
\newcommand\noot[1]
 {\marginpar{\small\begin{minipage}{25mm}\begin{flushleft}\thenootje : #1\end{flushleft}\end{minipage}}\addtocounter{nootje}{1}}
\begin{document}

\title{Asymptotics for partitions over the Fibonacci numbers and related sequences, II}


\author{Michael Coons}
\author{Simon Kristensen}
\author{Mathias L.~Laursen}

\address{Aarhus University, Ny Munkegade 118, 8000 Aarhus C, Denmark}
\email{mjc@math.au.dk, sik@math.au.dk, mathiaslaursen@hotmail.com}

\address{California State University, 400 West First Street, Chico, California 95929, USA}
\email{mjcoons@csuchico.edu}

%use when you want to break the address over multiple lines
%\newline
%\hspace*{\parindent}

\date{\today}

%%%%%%%%%%%%%%%%%%%%%%%%%%%%%%%%%%%%%%%%%%%%%%%%
\begin{abstract} In a recent paper of the same title, we considered nondistinct partitions of integers over the Fibonacci sequence and related linearly recurrent sequences. In this paper, continuing this story, we apply the same method, but consider distinct partitions. In particular, let $q_F(m)$ denote the number of ways of writing $m$ as the sum of distinct Fibonacci numbers. It is well-known that, while $q_F(m)$ behaves wildly, its partial sums behave in a more structured way. In particular, Kempton showed that the function $N^{-\log 2/\log\varphi}\sum_{m\leqslant N}q_F(m)$ is asymptotically logarithmically-periodic. In this paper, we improve on this result by providing an explicit error term. Additionally, we prove the analogous result for a class of more general linear recurrent sequences.
\end{abstract}
%%%%%%%%%%%%%%%%%%%%%%%%%%%%%%%%%%%%%%%%%%%%%%%%

\maketitle

%\vspace{-.7cm}
%%%%%%%%%%%%%%%%%%%%%%%%%%%%%%%%%%%%%%%%%%%%%%
\section{Introduction}
%%%%%%%%%%%%%%%%%%%%%%%%%%%%%%%%%%%%%%%%%%%%%%

Let $F_k$ denote the $k$th Fibonacci number, where, as usual, $F_0=0$ and $F_1=1$, and for $k\geqslant 1$, we have $F_{k+1}=F_k+F_{k-1}$. The object of this paper is the behaviour of the partial sums $\sum_{m\leqslant N}q_F(m)$, where $q_F(m)$ denotes the number of ways to write the integer $m$ as the sum of distinct Fibonacci numbers, that is, the number of solutions of \begin{equation}\label{eq:partitions} m=a_2F_2+a_3F_3+\cdots+a_kF_k+\cdots,\end{equation} where $a_k\in\{0,1\}$ for each $k\geqslant 2$. While the sequence $q_F(m)$ itself behaves erratically (see Chow and Slattery \cite{CS2021}), its partial sums are known to have determinable asymptotics. Indeed, Chow and Slattery showed that $\sum_{m\leqslant N}q_F(m)\asymp N^{\log 2/\log\varphi}$, where $\varphi=(1+\sqrt{5})/2$ is the golden mean. Then, using methods from dynamical systems, in an unpublished work, Kempton \cite{Kpre} showed that the function $N^{-\log 2/\log\varphi}\sum_{m\leqslant N}q_F(m)$ is asymptotically logarithmically-periodic. Recent work has also been done to understand the partial sums of the powers of $q_F(n)$. Indeed, Chow and Jones \cite{CJ2024} showed that $\sum_{m\leqslant N}q_F(m)^2\asymp N^{\log \lambda_2/\log\varphi}$ where $\lambda_2=2.48\cdots$ is the largest root of $X^3-2X^2-2X+2$, and, in general, Sanna \cite{S2025} showed for each positive integer $k$ that $\sum_{m\leqslant N}q_F(m)^k\asymp_k N^{\log \lambda_k/\log\varphi}$, where the $\lambda_k$ are effectively computable algebraic numbers, which satisfy $\lim_{k\to\infty}\lambda_k^{1/k}=\sqrt{\varphi}$.

In Part I of this work \cite{CKL2025}, using the saddle-point method combined with a result of Coons and Kirsten \cite{CK2009}, we determined asymptotics for the related quantity $p_F(n)$, the number of solutions to~\eqref{eq:partitions} where the condition $a_k\in\{0,1\}$ was replaced with the less strict condition $a_k\in\mathbb{N}\cup\{0\}$. Our method is dependent on the existence of invertible asymptotics for the related saddle point. Concerning this, in our paper, we noted that ``if one follows this method and tries to apply it to the distinct Fibonacci partition function $q_F(n)$ the resulting asymptotic is not monotonic---the leading term is oscillatory. This is precisely why this method doesn't immediately generalise to distinct Fibonacci partitions'' \cite[Remark 15]{CKL2025}. Well, it turns out, we were a little bit too quick with our judgement here. Indeed, the leading asymptotics contain an oscillatory term, but the oscillations are not so strong as to break monotonicity---colloquially, there's a bit of a wiggle, but just a little bit. With our new-found ability to overcome this wiggle, here, using number-theoretic methods, we give an explicit version of Kempton's result. 

\begin{theorem}\label{thm:mainpartitions} Let $q_F(m)$ be the number of partitions of $m$ over distinct Fibonacci numbers. Then, as $N\to\infty$, 
\vspace{-.2cm}
$$\sum_{m\leqslant N} q_F(m)=\psi_{F}\left(\frac{\log N}{\log\varphi}\right)N^{\frac{\log 2}{\log\varphi}}\left(1+O\left(\frac{1}{N}\right)\right),$$ where $\psi_{F}(x)$ is an explicit positive $1$-periodic function and $\varphi$ is the golden mean. 
\end{theorem}

Unlike in Kempton's description of $\psi_{F}$, we are able to determine $\psi_{F}$ explicitly, relating it to standard number theoretic functions, including the Euler $\Gamma$-function and the Riemann $\zeta$-function. 

\begin{figure}[h]
\begin{center}
  \includegraphics[width=0.975\linewidth]{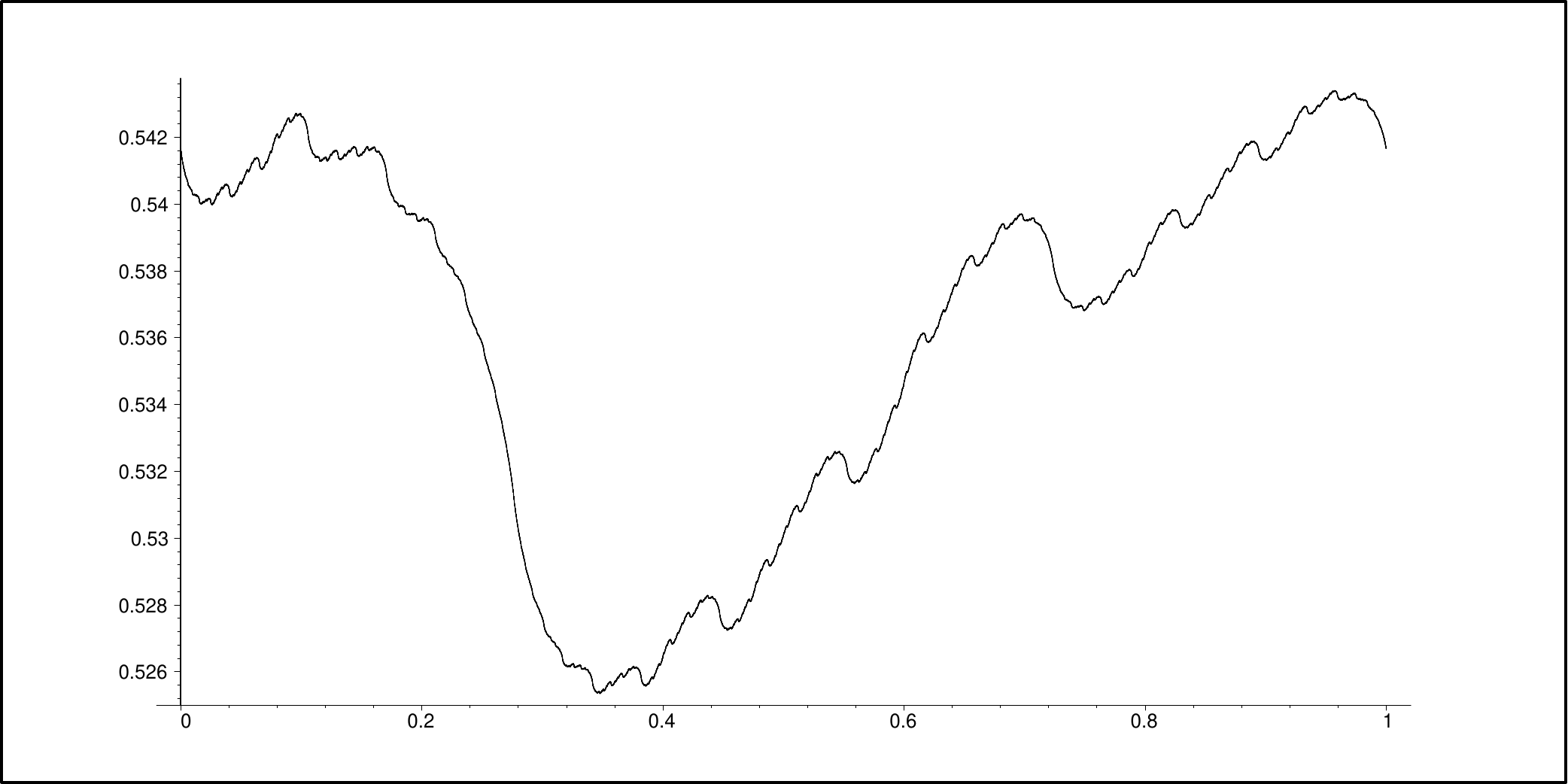}
\end{center}
\caption{The function $\psi_F(x)$ as  approximated by $N^{-\log 2/\log\varphi}\sum_{m\leqslant N}q_F(m)$ for $N\in[\varphi^{26},\varphi^{27}]$ and normalising to the unit interval $[0,1]$.}
\label{fig:1}
\end{figure}

Our method extends to cover linearly recurrent sequences akin to the Fibonacci sequence. In particular, suppose that $P_k$ is a strictly increasing linearly recurrent sequence of positive integers of degree at least $2$, such that the characteristic polynomial $\chi_P(x)$ of $P_k$ has a single dominant root $\beta>1$ and $P_k=\lambda\beta^k+\lambda_2\beta_2^k+\cdots+\lambda_r\beta_r^k,$ where $\lambda,\lambda_2,\ldots,\lambda_r$ are constants and $\beta_2,\ldots,\beta_r$ are the algebraic conjugates of $\beta$. 

\begin{theorem}\label{thm:genpart} Let $q_P(m)$ be the number of partitions of $m$ over distinct members of the sequence $P_k$. If $\beta\in(1,10]$, then, as $N\to\infty$, for any fixed $\varepsilon\in(0,1)$, $$\sum_{m\leqslant N} q_P(m)=\psi_{P}\left(\frac{\log N}{\log\beta}\right)N^{\frac{\log 2}{\log\beta}}\left(1+O\left(N^{-\min\{{\log(\beta/|\beta_2|)}/{\log\beta},1-\varepsilon\}}\right)\right),$$ where $\psi_{P}(x)$ is an explicitly determinable positive $1$-periodic function. 
\end{theorem}

\noindent We find the case of $\beta\in(1,2)$ to be the most interesting, since in that case $q_P(m)>1$ with natural density one. It's really just an artefact of our method that it does cover some $\beta>2$---in this case $q_p(m)=0$ with density one. 

Our paper is organised as follows. In the next section, we record the necessary asymptotic results to apply our method in the Fibonacci case. In Section \ref{sec:q}, we prove Theorem \ref{thm:mainpartitions}. In Section~\ref{sec:P}, we quickly discuss the proof of Theorem \ref{thm:genpart} as the argument follows the previous section nearly {\em mutatis mutandis}---the necessary results are recorded in this section. In our final section, we discuss further remarks and possible future directions.

\begin{remark}
Like with many asymptotic results, our arguments are inherently technical. Much of this technicality was discussed in depth in Part I, and so we have taken the liberty of assuming some familiarity with that paper. This paper can certainly be read on its own, but in combination with Part I, the picture is more complete. We make mention of Part I often, and have borrowed sentences and constructions from it freely. The nuance in this paper comes from understanding the nature of certain oscillations. Other than this, the arguments here are essentially those of Part I.
\end{remark}

%%%%%%%%%%%%%%%%%%%%%%%%%
%%%%%%%%%%%%%%%%%%%%%%%%%
\section{Necessary asymptotic results}
%%%%%%%%%%%%%%%%%%%%%%%%%
%%%%%%%%%%%%%%%%%%%%%%%%%

As stated in the Introduction, our proofs require several asymptotic results. In this section, for convenience, we collect the necessary results---these were also recorded or proven in Part I. 

A substantial part of our argument will be addressing asymptotics related to the product $s^{-z} \Gamma(z) \eta(1+z) \zeta_F(z)$, where $\Gamma$ is Euler's $\Gamma$-function, $\eta$ is Dirichlet's $\eta$-function, and $\zeta_F$ is the function defined on $\Re(z)>0$ by the  Dirichlet series $\sum_{k\geqslant 1}F_k^{-z}$. Note that here, $F_1=F_2=1$ is doubly represented---this will be important later. Also, $\eta(1+z)=(1-2^{-z})\zeta(1+z)$, where $\zeta$ is Riemann's $\zeta$-function. We will address the pieces of the above-mentioned product in order.

We start by noting that for $s > 0$, the function $z \mapsto s^{-z}$ can be expanded in a Taylor series
\begin{equation}
\label{eq:sto-z}
s^{-z} = 1 - z \log s + \frac{1}{2} (\log s)^2 z^2 + O(z^3),
\end{equation}
which converges for all $z \in \mathbb{C}$, and the $\Gamma$-function has simple poles at non-positive integers, where the residue at $-n$ is $$\res{s}{-n}{\Gamma(s)}=\frac{(-1)^n}{n!}.$$
The function $\eta(1+z)$ inherits the zeros and poles of the function $\zeta(1+z)$, as well as the extra zeros of $(1-2^{-z})$ which occur at $z=2\pi i n/\log2$ for all $n\in\mathbb{Z}$. The trivial zeros of $\eta(1+z)$ at all points in $-2\mathbb{N}-1$ will cancel the corresponding poles of $\Gamma$ in the product. Around $z=0$, we require the well-known asymptotics, 
\begin{equation}
\label{eq:gamma_0}
\Gamma(z) = \frac{1}{z} - \gamma + \frac{1}{2} \left(\gamma^2 + \frac{\pi^2}{6}\right)z + O(z^2),
\end{equation}
where $\gamma$ is the Euler--Mascheroni constant,
\begin{equation}
\label{eq:zeta_0}
\zeta(1+z) = \frac{1}z + \gamma - \gamma_1 z + O(z^2),
\end{equation}
where $\gamma_1$ is the first Stieltjes constant,
$\gamma_1 = \lim_{m \rightarrow \infty}  \big(\sum_{k=1}^m \frac{\log k}{k}  - \frac{1}{2} (\log m)^2\big)$, and \begin{equation}\label{eq:etaz=0}\eta(1+z)=\log 2 +\left(\gamma\log 2-\frac{(\log 2)^2}{2}\right)z+O(z^2).\end{equation} Equation \eqref{eq:gamma_0} can be proved by combining the functional equation for $\Gamma(z)$ with series 8.321.1 in Gradshteyn and Ryzhik \cite[p.~903]{GRtables}; see also Nielsen \cite[p.~40]{N1906}. The series \eqref{eq:zeta_0} is a classical result---the constant term $\gamma$ was given by Euler and the rest of the terms, known as the Stieltjes constants, were exactly determined by Stieltjes in a letter he wrote to Hermite in June 1885 \cite[p.151--155]{HS1905}. The series for $\eta(1+z)$ follows from its definition combined with \eqref{eq:zeta_0}.

Finally, we arrive at the slightly more complicated, or at least less well-known, function $\zeta_F$. Fortunately, this function was studied by Navas \cite{N2001} and its properties were further developed in Part I. We collect the necessary results from Navas's paper in the following proposition.

\begin{proposition}[Navas, 2001]\label{prop:navas1} The Dirichlet series $\zeta_F(z)$ can be continued analytically to a meromorphic function on all of $\mathbb{C}$, still called $\zeta_F(z)$, whose singularities are simple poles at $s=s(n,k)=-2k+\frac{\pi i(2n+k)}{\log\varphi},$ for $n,k\in\mathbb{Z}$ and $k\geqslant 0$, with $$\res{z}{s}{\zeta_F(z)}=\frac{(-1)^k5^{s/2}{-s\choose k}}{\log\varphi}.$$ Moreover, we have that $\zeta_F(-(4m+2))=0$ for each $m\in\mathbb{N}_0$, and $\zeta_F(-1)=-1$.
\end{proposition}

As in Part I, due to the product $s^{-z} \Gamma(z) \eta(1+z) \zeta_F(z)$ having a double pole at $z=0$, we require the following lemma \cite[Lemma 5]{CKL2025} taken from Part I.

\begin{lemma}\label{lem:zetaF_0} Near $z=0$, we have $$\zeta_F(z) = \frac{1}{\log\varphi}\cdot\frac{1}{z} + \frac{\log 5 - \log\varphi}{2\log\varphi} + \left(\frac{\log\varphi-3\log 5}{12} + \frac{(\log 5)^2}{8\log\varphi} + c_1 \right)z+O(z^2),$$ where $c_1:= \sum_{k\geqslant 1} \frac{(-1)^{k}}{k(\varphi^{2k} + (-1)^{k+1})}= -0.20436188\cdots.$
\end{lemma}

In addition to the special function results above, we require one further well-known asymptotic; in particular, that, as $s\to 0^+$, \begin{equation}\label{eq:e2s}-\log(1-e^{-2s})=-\log s-\log2+H(s),\end{equation} where $H(s)=O(s)$ is a power series---while it may seem strange not to just write $O(s)$ in this equations, at this point it is important to keep note that this function is power series, since later, when we differentiate, we will know that this term does not introduce any oscillations.

With these results recorded, we are now ready to proceed with the proof of our main result concerning distinct Fibonacci partitions.

%%%%%%%%%%%%%%%%%%%%%%%%%
%%%%%%%%%%%%%%%%%%%%%%%%%
\section{The normalised mean of distinct Fibonacci partitions}\label{sec:q}
%%%%%%%%%%%%%%%%%%%%%%%%%
%%%%%%%%%%%%%%%%%%%%%%%%%

In this section, we deal completely with the case of distinct partitions over Fibonacci numbers. Since the method of Part I is employed, we assume the reader is acquainted with that work and move a bit quicker here---ensuring we spend time on the differences. Towards this goal, recall that $q_F(m)$ is the number of ways to write the integer $m$ as the sum of distinct Fibonacci numbers (with only one 1). We define $$
G_2(x):=\sum_{N\geqslant 0}\left(\sum_{m\leqslant N} q_F(m)\right)x^N=(1-x)^{-1}\prod_{k\geqslant 2}\left(1+x^{F_k}\right)=(1-x^2)^{-1}\prod_{k\geqslant 1}\left(1+x^{F_k}\right)
$$
and consider the function $$G(e^{-s}):=\prod_{k\geqslant 1}\left(1+e^{-sF_k}\right).$$ Our main result will reduce to understanding the asymptotics of $\log G_2(e^{-s})$ as $s\to 0^+$, and applying the following result of Coons and Kirsten \cite{CK2009}, which originally appeared in Coons's master's thesis (written under Kirsten), and which itself was inspired by the work of Nanda \cite{N1954} and Richmond \cite{R1975,R1975/76}. 

\begin{theorem}[Coons and Kirsten, 2009]\label{thm:CK2009}
If $\Lambda(x)=\sum_{N\geqslant 0}\lambda(N)x^N$ is the generating function of a sequence $\lambda(N)$ that converges for $|x|<1$, and $s=\alpha$ is the solution of $\frac{d}{ds}\big(s+\frac{1}{N}\log \Lambda(e^{-s})\big)=0$, then, as $N$ tends to infinity,  $$\lambda(N)=\frac{e^{N\alpha}\Lambda(e^{-\alpha})}{\sqrt{2\pi}}\left(\sqrt{\frac{1}{-\left.\frac{dN}{ds}\right|_{s=\alpha}}}+O\left(\frac{1}{N^{3/2}}\right)\right).$$ Here, $\alpha$ must be thought of as being replaced by its large-$N$ asymptotic expansion so that the asymptotic of $\lambda(N)$ represents a large-$N$ asymptotic.
\end{theorem}

Building off of Mellin's formula for the exponential function, as we did with the function $F(e^{-s})$ in Part I, we have that \begin{equation}\label{eq:asymp2}\log G(e^{-s})=\frac{1}{2\pi i} \int_{a-i\infty}^{a+i\infty} s^{-z} \Gamma(z) \eta(1+z) \zeta_F(z) dz,
\end{equation}
where $a>0$, $\eta(z)=\sum_{m\geqslant 1}\frac{(-1)^{m+1}}{m^z}=(1-2^{1-z})\zeta(z)$ denotes the Dirichlet $\eta$-function. In Part I, the $\eta$-function's place was taken by the $\zeta$-function---here it is the `$+$'-sign in the products defining $G$ that causes the change. Note that $\eta(z)$ is entire and shares the trivial zeros of $\zeta(z)$. Considering values, we have that $\eta(1)=\log(2)$ and $\eta(0)=1/2$. Following the argument for $\log F(e^{-s})$ from Part I, the integral in \eqref{eq:asymp2} is equal to the the sum over the residues at the poles of the integrand.

\begin{theorem}\label{thm:Gform} The function $G(e^{-s})$ defined above satisfies, as $s\to0^+$, $$\log G(e^{-s})=(-\log s)\frac{\log 2}{\log\varphi}+\frac{\log 2\cdot \log(25/2\varphi^2)}{2\log\varphi}+\frac{s}{2}+\psi_0(s)+O(s^2),$$ where $\psi_0(s)$ satisfies $\psi_0(s)=\psi_0(\varphi s)$. 
\end{theorem}

\begin{proof} Combining our knowledge of the functions $s^{-z},$ $\Gamma(z)$ and $\eta(1+z)$ as recorded in the previous section, with Proposition~\ref{prop:navas1}, we observe that the integrand $s^{-z} \Gamma(z) \eta(1+z) \zeta_F(z)$ has 
\begin{itemize}
\item a double pole at $z=0$,
\item a simple pole coming from $\Gamma(z)$ at $z=-1$,
\item double poles at $z\in-4\mathbb{N}$, and
\item simple poles off the real line at $z=s(n,k)$, for $k\geqslant 0$ and $k\neq -2n$.
\end{itemize} where, as before the poles of $\Gamma(z)$ at $z=-(4m+2)$ for $m\in\mathbb{N}_0$ are cancelled by the corresponding zeros of $\zeta_F(z)$ and the poles of $\Gamma(z)$ at odd negative integers $z=-(2m+1)$ for $m\in\mathbb{N}$ are cancelled by the corresponding zeros of $\eta(1+z)$.

Applying \eqref{eq:sto-z}, \eqref{eq:gamma_0}, Lemma \ref{lem:zetaF_0} and \eqref{eq:etaz=0}, the double pole at $z=0$ contributes $$\res{z}{0}{s^{-z} \Gamma(z) \eta(1+z) \zeta_F(z)}=(-\log s)\frac{\log 2}{\log\varphi}+(2\log 5-2\log\varphi-\log 2)\frac{\log 2}{2\log\varphi}.$$ 

The simple pole at $z=-1$ contributes $$\res{z}{-1}{s^{-z} \Gamma(z) \eta(1+z) \zeta_F(z)}=s\, \eta(0)\, \zeta_F(-1)\cdot\res{z}{-1}{\Gamma(z)}=s\left(\frac{1}{2}\right)(-1)(-1)=\frac{s}{2}.$$ 

Using the proof of Prop.~8 from Part I, replacing the contribution $\zeta(1-4n)=B_{2n}/4n$ with the contribution  $\eta(1-4n)=(1-2^{4n})\zeta(1-4n)=(1-2^{4n})B_{2n}/4n$, we have
$$g_1(s):=\sum_{n=1}^\infty \res{z}{-4n}{s^{-z} \Gamma(z) \eta(1+z) \zeta_F(z)}=\sum_{n=1}^\infty \alpha_n\, s^{4n}-\log s\sum_{n=1}^\infty\beta_n s^{4n},$$ where $\beta_n:={b_{-4n}\eta(1-4n)}/{(4n)!}$, and $$\alpha_n:=\frac{(1-2^{4n})B_{2n}}{4n\cdot(4n)!}\left(b_{-4n}\left(\sum_{k=1}^{4n} \frac{1}{k} - \gamma\right) + c_{-4n} - \frac{5^{-2n}}{2}  \frac{(4n)!}{((2n)!)^2}\right)+\frac{b_{-4n}}{(4n)!}\cdot\eta'(1-4n),$$ where $B_{2n}$ is the $2n$-th Bernoulli number. Moreover, as $s\to0^+$, we have $g_1(s)=O(s^2)$. 

Finally, set
\begin{align*} h_1(s):=&\, \sum_{k=0}^\infty\sum_{n\in\mathbb{Z}\backslash\{-k/2\}} \res{z}{s(n,k)}{s^{-z} \Gamma(z) \eta(1+z) \zeta_F(z)}\\
=&\, \frac{1}{\log\varphi}\sum_{k=0}^\infty\sum_{n\in\mathbb{Z}\backslash\{-k/2\}} (-1)^k{-s(n,k)\choose k} \left(\frac{s}{\sqrt{5}}\right)^{-s(n,k)}\Gamma(s(n,k))\, \eta(1+s(n,k)).\end{align*} Now, since $|(1-2^{it})|<2,$ for all real $t$, we have, from the result for $\zeta(1+it)$ that $\eta(1+it)=O(\log |t|)$, so that $\psi_0(s):=[k=0]h_1(s)=O(1)$, where $[k=0]h_1(s)$ denotes the $k=0$ term of $h_1(s)$. For $k\neq 0$, we note that $$\big|\Gamma(s(n,k))\eta(1+s(n,k))\big|=\big|1-2^{-s(n,k)}\big|\cdot\big|\Gamma(s(n,k))\zeta(1+s(n,k))\big|<16^k\big|\Gamma(s(n,k))\zeta(1+s(n,k))\big|.$$ Thus, as in the proof of Prop.~9 from Part I, we obtain the analogue of the final displayed equation of that proof, with the added $16^k$ factor; that is, $$\big|\psi_0(s)\big|<\frac{d}{\big(e^{\frac{\pi^2}{2\log\varphi}}-1\big)}\left(\frac{16s^2}{10\pi\cdot\log\varphi\cdot \big(e^{\frac{\pi^2}{2\log\varphi}}-1\big)}\right)^{k},$$ where $d>0$ is independent of $n$ and $k$. 
This proves the first assertion. The second assertion, that $\psi_0(s)=\psi_0(\varphi s)$, follows exactly as in the proof in Prop.~9 of Part I.
\end{proof}

Combing Theorem \ref{thm:Gform} with \eqref{eq:e2s}, we have the following result.

\begin{corollary}\label{cor:Gform} The function $G_2(e^{-s})$ defined above satisfies,  as $s\to0^+$, $$\log G_2(e^{-s})=(-\log s)\left(\frac{\log 2}{\log\varphi}+1\right)+\frac{\log 2\cdot\log(25/2\varphi^2)}{2\log\varphi}-\log2+\psi_0(s)+H(s)+O(s^2).$$ \end{corollary}

At this point, we arrive at the main point of difference from Part I, the oscillatory function of interest, $\psi_0(s)$, where, as defined above,
$$\psi_0(s)=[k=0]h_1(s)=\frac{1}{\log\varphi}\sum_{n\in\mathbb{Z}\backslash\{0\}} \left(\frac{s}{\sqrt{5}}\right)^{-\frac{2\pi i n}{\log \varphi}}\Gamma\left(\frac{2\pi i n}{\log \varphi}\right)\, \eta\left(1+\frac{2\pi i n}{\log \varphi}\right).$$ 
It is clear that $\psi_0(s)=\psi_0(\varphi s)$ for all $s$. This function makes an important appearance in the equation defining the saddle point $s=\alpha$. In particular, the saddle point is the solution to $$\frac{d}{ds}\left(s+\frac{1}{N}\log G_2(e^{-s})\right)=0,$$ which, using the asymptotics in Corollary \ref{cor:Gform}, becomes \begin{align}
\label{eq:nG} N=-\left.\frac{d}{ds}\log G_2(e^{-s})\right|_{s=\alpha}&=\frac{1}{\alpha}\left(\frac{\log 2}{\log\varphi}+1\right)+\frac{\psi_1(\alpha)}{\alpha}-H'(0)+O(\alpha)\\ 
\label{eq:n}&=\frac{1}{\alpha}\left(\frac{\log 2}{\log\varphi}+1+\psi_1(\alpha)\right)\left(1+O(\alpha)\right),\end{align} where $$\psi_1(s):=s\frac{d}{ds}\psi_0(s)=\frac{1}{\log\varphi}\sum_{n\in\mathbb{Z}\backslash\{0\}} \left(\frac{s}{\sqrt{5}}\right)^{-\frac{2\pi i n}{\log \varphi}}\Gamma\left(1+\frac{2\pi i n}{\log \varphi}\right)\, \eta\left(1+\frac{2\pi i n}{\log \varphi}\right)$$ is an oscillatory function satisfying $\psi_1(s)=\psi_1(\varphi s)$. Now, our whole method relies on the ability to invert the relationship between $N$ and the saddle point $\alpha$. In particular, to apply the result of Coons and Kirsten, the quantity $\sqrt{-(dN/ds)^{-1}|_{s=\alpha}}$ must be well-defined. Thus, the asymptotics for $N$, as a function of $\alpha$, must be asymptotically strictly monotonic. To this end, noting that $$-\frac{d}{ds}\left(\frac{\psi_1(s)}{s}\right)=\frac{\psi_1(s)}{s^2}-\frac{\psi_1'(s)}{s}$$ and using \eqref{eq:nG}, we have \begin{equation}\label{eq:dnds}-\left.\frac{dN}{ds}\right|_{s=\alpha}=\frac{1}{\alpha^2}\left(\frac{\log 2}{\log\varphi}+1\right)+\frac{\psi_2(\alpha)}{\alpha^2}+O(1)=\frac{1}{\alpha^2}\left(\frac{\log 2}{\log\varphi}+1+\psi_2(\alpha)\right)\left(1+O(\alpha^2)\right),\end{equation} where $$\psi_2(s):=\psi_1(s)-s\psi_1'(s)=\sum_{n\in\mathbb{Z}\backslash\{0\}} \left(\frac{s}{\sqrt{5}}\right)^{-\frac{2\pi i n}{\log \varphi}}\left(\frac{1}{\log\varphi}+\frac{2\pi i n}{(\log\varphi)^2}\right)\Gamma\left(1+\frac{2\pi i n}{\log \varphi}\right)\, \eta\left(1+\frac{2\pi i n}{\log \varphi}\right)$$ is, as well, an oscillatory function satisfying $\psi_2(s)=\psi_2(\varphi s)$. But, even though $\psi_2(s)$ is oscillatory, it has very small oscillations. The following lemma implies that $\frac{\log 2}{\log\varphi}+1+\psi_2(s)>0$ for all $s$, and so the leading asymptotics of $N$, as a function of $s$ with $s\to 0^+$, are strictly monotonic, hence invertible.

\begin{lemma}\label{lem:psi2bound} For all $s>0$, we have $|\psi_2(s)|<0.0000046$.
\end{lemma}

\begin{proof} We have \begin{align*}|\psi_2(s)|&\leqslant \left(\frac{2\pi}{(\log\varphi)^2}+\frac{1}{\log\varphi}\right)\sum_{n\in\mathbb{Z}\backslash\{0\}} n\left|\Gamma\left(1+\frac{2\pi i n}{\log \varphi}\right)\, \eta\left(1+\frac{2\pi i n}{\log \varphi}\right)\right|\\
&\leqslant 2\left(\frac{2\pi}{(\log\varphi)^2}+\frac{1}{\log\varphi}\right)\frac{2\pi}{\log\varphi}\cdot 2\sum_{n=1}^\infty n^2\left|\Gamma\left(\frac{2\pi i n}{\log \varphi}\right)\, \zeta\left(1+\frac{2\pi i n}{\log \varphi}\right)\right|,\end{align*} since $|1-2^{-it}|<2$ for all real $t$ and $\Gamma(1+z)=z\Gamma(z)$ for all $z$. Note that for all nonzero integers $n$, $|2\pi i n/\log\varphi|>13$. In this range, an old result of Landau \cite{L1903} applies to give $$\left|\zeta\left(1+\frac{2\pi i n}{\log \varphi}\right)\right|\leqslant 2\log\left|\frac{2\pi n}{\log \varphi}\right|.$$ Also, for any real nonzero $t$, \begin{equation}\label{eq:absGamma}|\Gamma(it)|=\sqrt{\frac{\pi}{|t|\sinh(|t|)}}=\sqrt{\frac{2\pi}{|t|}}e^{-\pi |t|/2}\sqrt{\frac{1}{1-e^{-2\pi |t|}}},\end{equation} so that $$\left|\Gamma\left(\frac{2\pi i n}{\log \varphi}\right)\right|<\frac{7}{6}\sqrt{\frac{\log\varphi}{n}}\cdot e^{-\pi^2n/\log\varphi}.$$ Since $\log(2\pi n/\log\varphi)\cdot\sqrt{\log\varphi/n}<2$, for all $n\geqslant 1$, we have $$\left|\Gamma\left(\frac{2\pi i n}{\log \varphi}\right)\, \zeta\left(1+\frac{2\pi i n}{\log \varphi}\right)\right|<\frac{14}{3}\cdot e^{-\pi^2n/\log\varphi}.$$ Applying these to the above, we have that $$|\psi_2(s)|<\left(\frac{2\pi}{(\log\varphi)^2}+\frac{1}{\log\varphi}\right)\frac{2\pi}{\log\varphi}\cdot \frac{56}{3}\sum_{n=1}^\infty n^2\big(e^{-\pi^2/\log\varphi}\big)^n=0.00000456626214378943758\cdots,$$ where, to evaluate the series, we have used that $\sum_{n\geqslant 1}n^2x^n=x(1+x)/(1-x)^3$ is a rational function. This proves the result.
\end{proof}

This bound on $|\psi_2(s)|$ combined with \eqref{eq:n} and \eqref{eq:dnds} proves the following proposition.

\begin{proposition}\label{prop:alphadnda} For sufficiently small $\alpha>0$, or, equivalently, for sufficiently large $N$, $$\alpha=\frac{1}{N}\,\psi_3\left(\frac{\log N}{\log\varphi}\right)\left(1+O\left(\frac{1}{N}\right)\right),$$ and 
$$-\left.\frac{dN}{ds}\right|_{s=\alpha}=N^2\,\psi_4\left(\frac{\log N}{\log\varphi}\right)\psi_3\left(\frac{\log N}{\log\varphi}\right)^{-2}\left(1+O\left(\frac{1}{N^2}\right)\right),$$  where $$\psi_3(x):=\frac{\log 2}{\log\varphi}+1+\frac{1}{\log\varphi}\sum_{n\in\mathbb{Z}\backslash\{0\}} e^{-2\pi i n x}\,5^{\frac{\pi i n}{\log \varphi}}\,\Gamma\left(1+\frac{2\pi i n}{\log \varphi}\right)\, \eta\left(1+\frac{2\pi i n}{\log \varphi}\right)$$ and $$\psi_4(x):=\frac{\log 2}{\log\varphi}+1+\sum_{n\in\mathbb{Z}\backslash\{0\}} e^{-2\pi i n x}\,5^{\frac{\pi i n}{\log \varphi}}\,\left(\frac{2\pi i n}{(\log\varphi)^2}+\frac{1}{\log\varphi}\right)\Gamma\left(1+\frac{2\pi i n}{\log \varphi}\right)\, \eta\left(1+\frac{2\pi i n}{\log \varphi}\right)$$ are both $1$-periodic.
\end{proposition}

The final necessary piece, before applying Theorem~\ref{thm:CK2009}, is the following immediate corollary.

\begin{corollary}\label{cor:alphadnda} For sufficiently small $\alpha>0$, or equivalently, for sufficiently large $N$, $$-\log\alpha=\log N+\log\left(\psi_3\left(\frac{\log N}{\log\varphi}\right)\right)+O\left(\frac{1}{N}\right).$$
\end{corollary}

\begin{proof}[Proof of Theorem \ref{thm:mainpartitions}] As we have collected the requisite asymptotic results, we now evaluate the pieces of Theorem~\ref{thm:CK2009}, with the appropriate generating function and saddle points. Proposition~\ref{prop:alphadnda} gives both \begin{equation}\label{eq:S1}e^{N\alpha}=e^{\psi_3\left(\frac{\log N}{\log\varphi}\right)}\left(1+O\left(\frac{1}{N}\right)\right)\end{equation} and \begin{equation}\label{eq:S2}\left(\sqrt{\frac{1}{-\left.\frac{dN}{da}\right|_{a=\alpha}}}+O\left(\frac{1}{N^{3/2}}\right)\right)=\frac{1}{N}\psi_3\left(\frac{\log N}{\log\varphi}\right)\psi_4\left(\frac{\log N}{\log\varphi}\right)^{-1/2}\left(1+O\left(\frac{1}{N^{2}}\right)\right).\end{equation}

It remains to determine the small $\alpha$, and hence large $N$, asymptotics of $G_2(e^{-\alpha})$. To this end, we combine Corollary \ref{cor:Gform} with Corollary \ref{cor:alphadnda} to obtain, as $\alpha\to0^+$, or, equivalently, as $N\to\infty$, 
\begin{align*}\log G_2(e^{-\alpha})&=(\log N)\left(\frac{\log 2}{\log\varphi}+1\right)+\left(\frac{\log 2}{\log\varphi}+1\right)\log\left(\psi_3\left(\frac{\log N}{\log\varphi}\right)\right)\\
&\qquad\qquad+\frac{\log 2\cdot\log(25/2\varphi^2)}{2\log\varphi}-\log 2+\psi_5\left(\frac{\log N}{\log\varphi}\right)+O\left(\frac{1}{N}\right),
\end{align*} where $\psi_5(x)$ is the $1$-periodic function $$\psi_5(x):=\frac{1}{\log\varphi}\sum_{n\in\mathbb{Z}\backslash\{0\}} e^{-2\pi i nx}\,5^{\frac{\pi i n}{\log\varphi}}\,\Gamma\left(\frac{2\pi i n}{\log \varphi}\right)\, \eta\left(1+\frac{2\pi i n}{\log \varphi}\right),$$ determined by $\psi_0(s)$. Thus, we have \begin{equation*}G_2(e^{-\alpha})=\psi_6\left(\frac{\log N}{\log\varphi}\right)N^{\left(\frac{\log 2}{\log\varphi}+1\right)}\left(1+O\left(\frac{1}{N}\right)\right),\end{equation*} where \begin{equation*}\psi_6(x):=\exp\left(\left(\frac{\log 2}{\log\varphi}+1\right)\log\left(\psi_3\left(x\right)\right)+\frac{\log 2\cdot\log(25/2\varphi^2)}{2\log\varphi}-\log 2+\psi_5\left(x\right)\right).\end{equation*} Combining this with equations \eqref{eq:S1} and \eqref{eq:S2}, and setting \begin{equation}\label{eq:psiF} \psi_F(x):=\frac{e^{\psi_3(x)}\psi_3(x)\psi_6(x)}{\sqrt{2\pi\psi_4(x)}}\end{equation} finishes the proof of the theorem.
\end{proof}

%%%%%%%%%%%%%%%%%%%%%%%%%
%%%%%%%%%%%%%%%%%%%%%%%%%
\section{The normalised mean of distinct partitions over general linear recurrences with a dominant root}\label{sec:P}
%%%%%%%%%%%%%%%%%%%%%%%%%
%%%%%%%%%%%%%%%%%%%%%%%%%

As stated in the Introduction, we suppose that $P_k$ is a strictly increasing linearly recurrent sequence of positive integers of degree at least $2$, such that the characteristic polynomial $\chi_P(x)$ of $P_k$ has a single dominant root $\beta>1$ and $P_k=\lambda\beta^k+\lambda_2\beta_2^k+\cdots+\lambda_r\beta_r^k,$ where $\lambda,\lambda_2,\ldots,\lambda_r$ are constants and $\beta_2,\ldots,\beta_r$ are the algebraic conjugates of $\beta$. For such $P_k$, we wish to asymptotically understand the partial sums of the number $q_P(m)$ of solutions to \begin{equation}\label{eq:nP} m=a_1P_1+a_2P_2+\cdots+a_kP_K+\cdots,\end{equation} with $a_k\in\{0,1\}$. 

The case of partitions over $P_k$ is carried out exactly as in the case with the Fibonacci numbers $F_k$, but with the Fibonacci zeta function replaced by the zeta function $\zeta_P(z)$ which is the meromorphic continuation of the Dirichlet series $\sum_{k\geqslant 1}P_k^{-z}$ ($\Re(z)>0$). We consider the generating function $$G_P(x):=\sum_{N\geqslant 0}\left(\sum_{m\leqslant N}q_P(m)\right)x^N=(1-x)^{-1}\prod_{k\geqslant 1}\left(1+x^{P_k}\right).$$
Continuing the analogy with the Fibonacci partitions, we need asymptotic results, as $s\to0^+$, of the function $$\log G_P(e^{-s}) = \frac{1}{2\pi i} \int_{a-i\infty}^{a+i\infty} s^{-z} \Gamma(z) \eta(1+z) \zeta_P(z) dz.$$ Here, instead of the result of Navas used previously, we appeal to the following generalisation of Serrano Holdago and Navas Vicente~\cite{HN2023}.

\begin{proposition}[Serrano Holdago and Navas Vicente, 2023]\label{prop:HN2023} Let $P(x)$ be the minimal polynomial with $\deg P(x)=r$ of the linear recurrence $P_k$ (as described above) and let $\beta>1$ be the dominant root of $P(x)$. Then the Dirichlet series $\zeta_P(z):=\sum_{k\geqslant 1}P_k^{-z}$ ($\Re(z)>0$) can be analytically continued to a meromorphic function, also denoted $\zeta_P(z)$, all of whose singularities are simple poles at the points $$s{(n,\boldsymbol{k})}:=\frac{\log|\beta^{-k_1}\beta_2^{k_1-k_2}\cdots \beta_r^{k_{r-1}}|}{\log \beta}+i\cdot\frac{\arg(\beta_1^{-k_1}\beta_2^{k_1-k_2}\cdots \beta_r^{k_{r-1}})+2\pi n}{\log \beta},$$ where $\boldsymbol{k}=(k_1,\ldots,k_{r-1})$, $\beta_2,\ldots,\beta_r$ are the algebraic conjugates of $\beta$, $n\in\mathbb{Z}$ and the parameters $k_1,\ldots,k_{r-1}$ are integers satisfying $0\leqslant k_{r-1}\leqslant k_{r-2}\leqslant \cdots \leqslant k_{1}$. 
\end{proposition}

As well, in the same paper, Serrano Holdago and Navas Vicente provided the values of the residues. The important one for our context is $$\res{z}{s(n,\boldsymbol{0})}{\zeta_P(z)}=\frac{\lambda^{s(n,\boldsymbol{0})}}{\log\beta}.$$ 

Replacing $\beta$ with $\varphi$ and $\lambda$ with $1/\sqrt{5}$, the arguments of the previous section go through immediately and with only trivial changes, but for one limitation:~we need of an analogue of Lemma~\ref{lem:psi2bound}! In particular, we must show that the function $\psi_2(s)$ associated to $P$, which we denote $\psi_{2,P}(s)$, satisfies $$|\psi_{2,P}(s)|<\frac{\log 2}{\log\beta}+1$$ for all $s$. The main obstacle to producing a completely general result using our method is the requirement of an explicit bound on $|\zeta(1+it)|$, with $t=2\pi n/\log\beta$. As $\beta$ increases, our argument necessitates a bound that works close to $t=0$, and, using the known results, say with $|t|>e$, this puts an upper bound on $\beta$.

\begin{lemma}\label{lem:psi2Pbound} For each $\beta\in(1,10]$ and all $s>0$, we have $|\psi_{2,P}(s)|<0.305310893739$.
\end{lemma}

\begin{proof} Following the proof for $\psi_2(s)=\psi_{2,F}(s)$, we have \begin{equation}\label{eq:psigammazeta}|\psi_{2,P}(s)|\leqslant 2\left(\frac{2\pi}{(\log\beta)^2}+\frac{1}{\log\beta}\right)\frac{2\pi}{\log\beta}\cdot 2\sum_{n=1}^\infty n^2\left|\Gamma\left(\frac{2\pi i n}{\log \beta}\right)\, \zeta\left(1+\frac{2\pi i n}{\log \beta}\right)\right|.\end{equation} We use a stricter, and much newer, result on the bounds of $\zeta(1+it)$ for real $t$. Indeed, recently, Hoo and Teo~\cite{HT2026}, have shown that for $|t|>e$, we have $|\zeta(1+it)|\leqslant 0.6443\log t.$ Thus, taking $\beta\leqslant 10$ to ensure that $|2\pi i n/\log\beta|>e$, this, along with \eqref{eq:absGamma}, gives $$\left|\Gamma\left(\frac{2\pi i n}{\log \beta}\right)\, \zeta\left(1+\frac{2\pi i n}{\log \beta}\right)\right|< 1.188265772\cdot e^{-\pi^2n/\log\beta},$$ since for $\beta\in(1,10]$ and $n\geqslant 1$, we have both $$\sqrt{1/(1-\exp(-4\pi^2n/\log\beta))}< 1.000000017902\quad\mbox{and}\quad\sqrt{\frac{\log\beta}{n}}\cdot \log\left(\frac{2\pi n}{\log\beta}\right)<1.84427402,$$ where the first inequality is closest when $n=1$ and $\beta=10$, and the second is closest when $n=1$ and $\beta=2.34043465\cdots$. Applying these as before, we have $$|\psi_{2,P}(s)|<\left(\frac{2\pi}{(\log\beta)^2}+\frac{1}{\log\beta}\right)\frac{8\pi}{\log\beta} (1.188265772)\sum_{n=1}^\infty n^2\big(e^{-\pi^2/\log\beta}\big)^n<0.3053108937392\cdots,$$ where the upper bound is closest when $\beta=10$. This proves the result.
\end{proof}

Using the bound in Lemma \ref{lem:psi2Pbound}, Theorem \ref{thm:genpart} follows {\em mutatis mutandis} from the previous section.

%%%%%%%%%%%%%%%%%%%%%%%%%
%%%%%%%%%%%%%%%%%%%%%%%%%
\section{Further remarks}\label{sec:C}
%%%%%%%%%%%%%%%%%%%%%%%%%
%%%%%%%%%%%%%%%%%%%%%%%%%

In this paper, we showed that for certain linearly recurrent sequences, the partial sums of the number of distinct partitions $\sum_{m\leqslant N}q_P(m)$ over that sequence, satisfy an asymptotic of the form \begin{equation}\label{eq:ass}\psi\left(\frac{\log N}{\log \beta}\right)N^{\frac{\log 2}{\log \beta}}(1+o(1)),\end{equation} where $\beta\in(1,10]$ is the unique dominant eigenvalue associated to the sequence $P_k$ and $\psi(x)$ is a (possibly) oscillatory $1$-periodic function. Of course, when $\beta=2$, $\psi(x)=1$ for all $x$---hence our use of `(possibly)' here. The most interesting part of our result is when $\beta\in(1,2)$ since then we have $q_P(m)>1$ most of the time. To extend our result, one could try to extend the range of $\beta$. Actually, the upper bound of $10$ for $\beta$ in Theorem \ref{thm:genpart} is quite close to what Hoo's and Teo's~\cite{HT2026} result on $\zeta$ can handle. The condition that $|2\pi i/\log\beta|>e$ is equivalent to $\beta<e^{2\pi/e}=10.089\cdots.$ So, one would need some other way of bounding $\psi_{2,P}(s)$ to continue using our method. One possibility is splitting the sum in \eqref{eq:psigammazeta} and considering those terms $n$ with $2\pi n/\log\beta<e$ separately; for those, one could try to apply the standard expansion of $\zeta(z)$ around $z=1$ and see what happens. Of course, it is immediate that for any $\beta>1$, we have $\sum_{m\leqslant N}q_P(m)\asymp N^{\log 2/\log\beta}$, since we are choosing from the two letters $\{0,1\}$ and asking (essentially) for a $\beta$-expansion---we have $P_k\sim\lambda\beta^k$. The problem is then understanding the function $\psi(x)$. 

Asymptotics of the form \eqref{eq:ass} are somewhat expected when one is dealing with the partial sums of regular sequences---sequences that are generalisations of linearly recurrent sequences. Recall that a linearly recurrent sequence is determined by the powers of one matrix; for example, for the Fibonacci sequence, $$F_k=(1\ 0)\begin{pmatrix}1&1\\ 1&0\end{pmatrix}^k\begin{pmatrix} 0\\ 1\end{pmatrix}.$$ A sequence $f$ is {\em regular} if there exist integers $k\geqslant 2$ and $d\geqslant 1$, vectors ${\bf u},{\bf v}\in\mathbb{C}^{d\times 1}$ and matrices ${\bf A}_0,\ldots,{\bf A}_{k-1}\in\mathbb{C}^{d\times d}$ such that $$f(n)={\bf u}^T {\bf A}_{i_0}{\bf A}_{i_1}\cdots {\bf A}_{i_s}{\bf v},$$ where $(n)_k=i_s\cdots i_1i_0$ is the base-$k$ expansion of $n$. It is quite common for the partial sums of regular sequences to satisfy an asymptotic of the form \eqref{eq:ass}; for examples, see Dumas \cite{D2013,D2014} and Coons and Lind \cite{CL2024}. As an interesting case, it may be worth trying to use our method to understand the behaviour of the ($\ell$-regular) partition functions $b_{k,\ell}(n)$, the number of partitions of the integer $n$ into powers of $\ell$, where each part appears at most $k$ times. In particular, can our method be used to address the conjectures of Calkin and Wilf \cite{CWpre} concerning $b_{k,\ell}(n)$? Here, one may have to consider a Dirichlet series $\sum_{n\geqslant 1}f(n)^{-z}$ where $f(n)$ is a $k$-regular sequence---certainly an object of interest in its own right.

%%%%%%%%%%%%%%%%%%%%%%%%%%%%%%%%%%%%%%%%%%%%%%
\subsection*{Acknowledgements}
%%%%%%%%%%%%%%%%%%%%%%%%%%%%%%%%%%%%%%%%%%%%%%
The work of M.~Coons was supported by Villum Fonden research grant VIL83378 and by the David W.~and Helen E.~F.~Lantis Endowment.

%%%%%%%%%%%%%%%%%%%%%%%%%%%%%%%%%%%%%%%%%%%%%%%%
\bibliographystyle{amsplain}

\begin{thebibliography}{99}

%\bibitem{dB1948}
%N.~G. de~Bruijn, \emph{On {M}ahler's partition problem}, Indagationes Math.~\textbf{10} (1948), 210--220.

\bibitem{HS1905}
B.~Baillaud (ed.) and H.~Bourget (ed.),
Correspondance d'Hermite et de Stieltjes, Tome I (8 novembre 1882 – 22 juillet 1889), Gauthier-Villars. xxii+ 477 p., 1905.

\bibitem{CWpre}
N.~J.~Calkin and H.~S.~Wilf, \emph{Binary partitions of integers and Stern-Brocot-like trees}, unpublished (1998), updated version August 5, 2009, 19 pages.

\bibitem{CS2021}
S.~Chow and T.~Slattery, \emph{On Fibonacci partitions}, J. Number Theory \textbf{225} (2021), 310--326.

\bibitem{CJ2024}
S.~Chow and O.~Jones, \emph{On the variance of the Fibonacci partition function}, J. Number Theory \textbf{257} (2024), 341--353.

\bibitem{CK2009}
M.~Coons and K.~Kirsten, \emph{General moment theorems for nondistinct unrestricted partitions}, J.~Math.~Phys.~\textbf{50} (2009), Article 013517.

\bibitem{CKL2025}
M.~Coons, S.~Kristensen, and M.~L.~Laursen, \emph{Asymptotics for partitions over the Fibonacci numbers and related sequences}, Comb.~Number Theory \textbf{14} (2025), no.~2, 119--140.

\bibitem{CL2024}
M.~Coons and J.~Lind, \emph{Radial asymptotics of generating functions of $k$-regular sequences}, Bull. Aust. Math. Soc. \textbf{110} (2024), no.~3, 528--534.

\bibitem{D2013}
Ph.~Dumas, \emph{Joint spectral radius, dilation equations, and asymptotic behavior of radix-rational sequences}, Linear Algebra Appl.~\textbf{438} (2013), no.~5, 2107--2126.

\bibitem{D2014}
Ph.~Dumas, \emph{Asymptotic expansions for linear homogeneous divide-and-conquer recurrences: algebraic and analytic approaches collated}, Theoret. Comput. Sci. \textbf{548} (2014), 25--53.
  
\bibitem{GRtables}
I.~S.~Gradshteyn and I.~M.~Ryzhik,
Table of integrals, series, and products, edited by D.~Zwillinger and V.~Moll, 8th Edition, Elsevier/Academic Press, xlv + 1133 p., 2015.

%\bibitem{HR1918}
%G.~H.~Hardy and S.~Ramanujan, \emph{Asymptotic formul\ae\ in combinatorial analysis}, Proc. Lond. Math. Soc. (2) \textbf{17} (1918), 75--115.

\bibitem{HT2026} 
E.~Q.~Hoo and L.-P.~Teo, \emph{Explicit bounds of $|\zeta(1+it)|$}, Rend. Semin. Mat. Univ. Padova \textbf{156} (2026), 201--233.

\bibitem{Kpre}
T.~Kempton, \emph{The dynamics of the Fibonacci partition function}, unpublished (2023). {\tt arXiv:2311.06006}

\bibitem{L1903}
E.~Landau, \emph{Neuer Beweis des Primzahlsatzes und Beweis des Primidealsatzes}, Math.~Ann.~\textbf{56} (1903), no.~4, 645--670.

%\bibitem{M1940}
%K.~Mahler, \emph{On a special functional equation}, J.~London Math.~Soc.~\textbf{15} (1940), 115--123. 

\bibitem{N1954}
V.~S.~Nanda, \emph{Partition theory and thermodynamics
of multidimensional oscillator assemblies}, Proc. Camb. Phil. Soc. \textbf{47} (1954), 591--601.

\bibitem{N2001}
L.~Navas, \emph{Analytic continuation of the {Fibonacci} {Dirichlet} series},
  Fibonacci Q. \textbf{39} (2001), no.~5, 409--418.

\bibitem{N1906}
N.~Nielsen, Handbuch der Theorie der Gammafunktion, Teubner, Leipzig, 1906.
  
%\bibitem{O1974} 
%F.~W.~J.~Olver, {Asymptotics and Special Functions}, Academic Press, New York, 1974.

\bibitem{R1975}
L.~B.~Richmond, \emph{The moments of partitions I},
Acta Arith.~\textbf{26} (1975), 411--425.

\bibitem{R1975/76}
L.~B.~Richmond, \emph{The moments of partitions II},
Acta Arith.~\textbf{28} (1975/76), 229--243.

\bibitem{S2025}
C.~Sanna, \emph{A note on the power sums of the number of Fibonacci partitions}, Discrete Anal. (2025), Paper No. 2, 13 pages.

\bibitem{HN2023}
\'{A}.~Serrano Holgado and L.~M.~Navas Vicente, \emph{The zeta function of a recurrence sequence of arbitrary degree}, Mediterr. J. Math. \textbf{20} (2023), no.~4, Paper No.~224, 18 pp.

%\bibitem{T1986}
%E.~C.~Titchmarsh, The theory of the Riemann zeta-function, 2 ed., Oxford University Press, New York, 1986.

\end{thebibliography}

%%%%%%%%%%%%%%%%%%%%%%%%%%%%%%%%%%%%%%%%%%%%%%%%

\end{document}